\documentclass[reqno,12pt]{amsart}
\usepackage{geometry}
\usepackage{epsfig}
\usepackage{amsmath}
\usepackage{amsfonts}
\usepackage{mathrsfs}
\usepackage{cite}
\usepackage[colorlinks,
            linkcolor=blue,
            anchorcolor=blue,
            citecolor=blue
            ]{hyperref}

\numberwithin{equation}{section} 
\allowdisplaybreaks[4]

\begin{document}
\newtheorem{theorem}{Theorem}[section]
\newtheorem{lemma}[theorem]{Lemma}
\newtheorem{corollary}[theorem]{Corollary}
\newtheorem{conjecture}[theorem]{Conjecture}
\newtheorem{question}[theorem]{Question}
\newtheorem{pro}[theorem]{Proposition}
\newtheorem{remark}[theorem]{Remark}
\newtheorem{open problem}[theorem]{Open problem}
\newcommand{\arcsinh}{\mbox{arcsinh}}
\newcommand{\arccosh}{\mbox{arccosh}}
\renewcommand{\theequation}{\arabic{section}.\arabic{equation}}
\newcommand{\IN}{\mathbb{N}}
\newcommand{\IR}{\mathbb{R}}
\newcommand{\IC}{\mathbb{C}}
\newcommand{\e}{\mbox{e}}
\newcommand{\K}{\mathscr{K}}
\newcommand{\E}{\mathscr{E}}
\newcommand{\La}{\mathscr{L}}
\newcommand{\s}{\sum_{n=0}^{\infty}}
\newcommand{\su}{\sum_{n=1}^{\infty}}
\newcommand{\la}{\lambda}
\newcommand{\G}{\Gamma}
\newcommand{\g}{\gamma}
\newcommand{\p}{\psi}

\def\square{\hfill${\vcenter{\vbox{\hrule height.4pt \hbox{\vrule width.4pt
height7pt \kern7pt \vrule width.4pt} \hrule height.4pt}}}$}

\title[]
{Monotonicity properties of  Gaussian hypergeometric functions with respect to the Parameter$^*$}

\author{Qi Bao$^{1}$}

\address{Qi Bao$^{1}$, $^{1}$ School of Mathematical Sciences, East China Normal University, Shanghai, 200241, China}

\email{52205500010@stu.ecnu.edu.cn}

\author{Miao-Kun Wang$^{2,**}$}

\address{Miao-Kun Wang$^{2}$, $^{2}$Department of Mathematics, Huzhou University, Huzhou, Zhejiang, 313000, China}

\email{wmk000@126.com; wangmiaokun@zjhu.edu.cn}

\author{Song-Liang Qiu$^{3,4}$}

\address{Song-Liang Qiu$^{3,4}$, $^{3}$ Department of Mathematics, Zhejiang Sci-Tech University, Hangzhou 310018, China; $^{4}$ Department of Mathematics, Lishui University, Lishui 323000, Zhejiang, China}
\email{sl$\_$qiu@zstu.edu.cn}



\thanks{$^*$The research was supported by the Natural Science Foundation of China(11701176, 11971142)}

\thanks{$^{**}$Corresonding author. Email: wmk000@126.com; wangmiaokun@zjhu.edu.cn. Tel.: +86 572 2321510; Fax: +86 572 2321163}

\begin{abstract}
The authors establish the necessary and sufficient conditions under which certain combinations of Gaussian hypergeometric function and elementary function are monotone in the parameter, which generalize the recent results of  generalized elliptic integrals of the first and second kinds obtained by Qiu et al. Moreover, the authors also prove two monotonicity theorems of generalized elliptic integrals from another point of view.
\end{abstract}
\maketitle

\section{introduction}
\setcounter{equation}{0}

Throughout this paper, we always let $r^{\thinspace\prime}=\sqrt{1-r^2}$ for $r\in [0,1]$, denote by $\IN$ ($\IR$) the set of positive integers (real numbers, respectively), and set $\IN_0=\IN\cup\{0\}$. For complex number $x$ with ${\rm Re}\,x>0$, let
\begin{align}\label{GAMMA1}
\Gamma(x)=\int_0^{\infty}t^{x-1}e^{-t}dt, \,\,
B(x,y)=\frac{\G(x)\G(y)}{\G(x+y)}, \,\,
\psi(x)=\frac{\Gamma'(x)}{\Gamma(x)}
\end{align}
be the classical Euler gamma, beta and psi (digamma) functions, respectively (cf. \cite{AS,AAR,Ask2,AVV0}). For complex numbers $a, b$ and $c$ with $c\neq 0,-1,-2,\cdots$, the Gaussian hypergeometric function is defined by
\begin{align}\label{2F1}
F(a,b;c;x)={}_2F_1(a,b;c;x)=\sum_{n=0}^{\infty}\frac{(a,n)(b,n)}{(c,n)}\frac{x^n}{n!} ~(|x|<1),
\end{align}
where $(a,n)$ is the Pochhammer symbol or shifted factorial defined as $(a,0)=1$ for $a\neq 0$, and $(a,n)=a(a+1)(a+2)\cdots (a+n-1)=\G(n+a)/\G(a)$ for $n\in\IN$ (cf. \cite{AS,AAR,Ask2,Be2,OLB}). $F(a,b;c;x)$ is said to be zero-balanced if $c=a+b$, and it converges absolutely for all $|x|<1$\, (cf. \cite[Theorem 2.1.1]{AAR}). It is well known that $F(a,b;c;x)$ has many important applications in several branches of mathematics, physics and engineering, and many other special functions in mathematical physics and even some elementary functions are particular or limiting cases of this function.

In 1769, Euler gave an important integral representation of Gaussian hypergeometric function ${}_2F_1$ (cf. \cite[15.3.1]{AS}). If ${\rm Re}\,c>{\rm Re} \,b>0$ and $x\in\IC\backslash[1,+\infty)$, then
\begin{align}\label{int}
{}_2F_1(a,b;c;x)=\frac{2\G(c)}{\G(b)\G(c-b)}
\int_0^{\pi/2} (\sin u)^{2b-1}(\cos u)^{2c-2b-1}(1-x\sin^2u)^{-a}du.
\end{align}

For $a\in (0,1)$ and $r\in (0,1)$, the generalized elliptic integrals of the first and second kinds are defined as
\begin{align}
\begin{cases}\label{Ka}
\K_a=\K_a(r)={\pi}F\left(a,1-a;1;r^2\right)/2,\\
\K_a'=\K_a'(r)=\K_a(r^{\thinspace\prime}),\\
\K_a(0)=\pi/2,\K_a(1)=\infty,
\end{cases}
\end{align}
and
\begin{align}
\begin{cases}\label{Ea}
\E_a=\E_a(r)={\pi}F\left(a-1,1-a;1;r^2\right)/2,\\
\E_a'=\E_a'(r)=\E_a(r^{\thinspace\prime}),\\
\E_a(0)=\pi/2,\E_a(1)=[\sin(\pi a)]/[2(1-a)],
\end{cases}
\end{align}
respectively (cf. \cite{AVVb,AQVV,HVV}). Taking $a=1/2$, $\K_{1/2}\equiv\K$ and $\E_{1/2}\equiv\E$ are the complete elliptic integrals of the first and second kinds, respectively (cf. \cite[17.3.9--17.3.10]{AS}).

During the past few years, both complete elliptic integrals ($\K$ and $\E$) and generalized elliptic integrals ($\K_a$ and $\E_a$) have been widely studied and applied in the theories of conformal invariants, quasiconformal mappings and Ramanujan's modular equations\cite{AQVV,Wang2,WCQ7,AlzerQiu1,HQM2019,YangC,WHC2020,ZT H1,yzh1,WangQC1,QMB,AB1,ZXH2,Zhang2,Qiu6,Qiu3,QMC1,G4,G5,Yang-Tian-RJ-2019,WCZ-2019-RJ,WLC-2018-RJ,MQJ-2021-JMAA}.
One of the meaningful tasks is to investigate the dependence on the parameter $a$ in $\K_{a}$ and $\E_{a}$, and thus show some stabilities of $\K_{a}$ and $\E_{a}$ with respect to $a$ and establish several sharp bounds for $\K_{a}$ and $\E_{a}$ in terms of $\K$ and $\E$.

In 2000, Anderson, Qiu, Vamanamurthy and Vuorinen \cite[Theorem 7.2]{AQVV} proved the following Theorem {\ref{theorem1}}.
\begin{theorem}\label{theorem1}\rm{}
For each $r\in(0,1)$, let $f,\,g$ be defined on $[0,1]$ by $\mu(a)=F(a-1,1-a;1;r^2)$ and $\nu(a)=F(a,1-a;1;r^2)$.

(1) If $1/2\leq a<b\leq1$, then all coefficients are positive in the Taylor series for $\mu(b)-\mu(a)$ in powers of $r^2$.

(2) If $0\leq a<b\leq1-1/\sqrt2$, then all coefficients are negative in the Taylor series for $\mu(b)-\mu(a)-(b-a)(2-a-b)r^2$ in powers of $r^2$.

(3) If $0\leq a<b\leq1/2$ $( 1/2\leq a<b\leq1, \mbox{respectively})$, then all coefficients are positive (negative, respectively) in the Taylor series for $\nu(b)-\nu(a)$ in powers of $r^2$.
\end{theorem}

It is apparent from Theorem \ref{theorem1} that $\K_a$ and $\E_a$ are both strictly increasing with respect to the parameter $a\in(0,1/2]$. Recently, Qiu, Ma and Bao\cite[Theorems 1.1 and 1.2]{QMB} presented the necessary and sufficient conditions under which certain familiar combinations, defined in terms of $\K_a$, $\E_a$ and elementary functions, are monotone in $a\in(0,1/2]$, so that some known related results were proved substantially. For example, they proved the following monotonicity theorem.
\begin{theorem}\label{exampleth1} \rm{}
Let $1.118763390276<\la_0<1.118763390286$ be as in \cite[Lemma 2.1\,(3)]{QMB}, for each $r\in(0,1)$ and $\la\in \IR$, define the functions $\varphi_{1,\,\la}$, $\varphi_{2,\,\la}$, $\varphi_{3,\,\la}$ and $\varphi_{4,\,\la}$  on $(0,1/2]$ by
\begin{align*}
\varphi_{1,\,\la}(a)&=\frac{\K_a(r)-\pi/2}{a^{\la}},\quad
\varphi_{2,\,\la}(a)=\frac{\pi/2-\E_a(r)}{a^{\la}}, \\
\varphi_{3,\,\la}(a)&=\frac{\K_a(r)-\E_a(r)}{a^{\la}} \quad \mbox{and} \quad
\varphi_{4,\,\la}(a)=\frac{\E_a(r)-r^{\thinspace\prime 2}\K_a(r)}{a^{\la}},
\end{align*}
respectively. Then the following conclusions hold:

(1) $\varphi_{1,\,\la}$ is strictly increasing (decreasing) on $(0,1/2]$ if and only if $\la \leq 0$ ($\la \geq 1$, respectively).

(2) $\varphi_{2,\,\la}$ is strictly increasing (decreasing) on $(0,1/2]$ if and only if $\la \leq -2$ ($\la \geq 0$, respectively).

(3) $\varphi_{3,\,\la}$ is strictly increasing (decreasing) on $(0,1/2]$ if and only if $\la\leq-1$ ($\la\geq1$, respectively).

(4) $\varphi_{4,\,\la}$ is strictly increasing (decreasing) on $(0,1/2]$ if and only if $\la\leq 1$ ($\la\geq\la_0$, respectively).
\end{theorem}

Besides, the authors \cite[Theorem 6.2]{QMB} further proved

\begin{theorem}\label{exampleth2}\rm{} For each $\la\in\IR$, $n\in\IN_0$ and for each $r\in(0,1)$, define the functions $\varphi_{5,\,\la}$ and $\varphi_{6,\,\la}$ on $(0,1/2]$ by
\begin{align*}
\varphi_{5,\,\la}(a)=\frac{\K_a(r)-\frac{\pi}{2}\sum_{k=0}^n \frac{(a,k)(1-a,k)}{(k!)^2}r^{2k}}{a^{\la}}
\quad \text{and}\quad
\varphi_{6,\,\la}(a)=\frac{\frac{\pi}{2}\sum_{k=0}^n \frac{(a-1,k)(1-a,k)}{(k!)^2}r^{2k}-\E_a(r)}{a^{\la}},
\end{align*}
respectively. Then the following conclusions hold:

(1) $\varphi_{5,\,\la}$ is strictly increasing (decreasing) on $(0,1/2]$ if and only if $\la\leq0$ ($\la\geq 1$, respectively).

(2) $\varphi_{6,\,\la}$ is strictly increasing (decreasing) on $(0,1/2]$ if and only if $\la\leq -1-\frac{1}{2n+1}$ ($\la\geq 0$ if $n=0$ or $\la\geq 1$ if $n\in \IN$, respectively).
\end{theorem}

Substituting $n=0$ in Theorem \ref{exampleth2},  Theorem \ref{exampleth2} (1) and (2) reduce to Theorem \ref{exampleth1} (1) and (2), respectively. In this paper, we shall consider more general situations. For fixed $c\in(0,\infty)$, let $a\in(0,c/2]$, $\la\in\IR$ and $x\in (0,1)$, then we define
\begin{equation}\label{fla12}
f_{1,\,\la}(a)=\frac{F(a,c-a;c;x)-1}{a^{\la}}, \quad f_{2,\,\la}(a)=\frac{1-F(a-1,c-a;c;x)}{a^{\la}},
\end{equation}
\begin{equation}\label{fla3}
f_{3,\,\la}(a)=\frac{F(a,c-a;c;x)-F(a-1,c-a;c;x)}{a^{\la}}
\end{equation}
and
\begin{equation}\label{fla4}
f_{4,\,\la}(a)=\frac{F(a-1,c-a;c;x)-(1-x)F(a,c-a;c;x)}{a^{\la}}.
\end{equation}
Obviously, when $c=1$ and $x=r^2$, $f_{1,\,\la}$, $f_{2,\,\la}$, $f_{3,\,\la}$ and $f_{4,\,\la}$ become $\varphi_{1,\,\la}$, $\varphi_{2,\,\la}$, $\varphi_{3,\,\la}$ and $\varphi_{4,\,\la}$ in Theorem \ref{exampleth1}, respectively. Naturally, the following Question \ref{question} is proposed.

\begin{question}\label{question}
For what values of $\la\in\IR$, $f_{i,\,\la}\,\, (i=1,2,3,4)$ are increasing (or decreasing) on $a\in(0,c/2]$?
\end{question}

Motivated by the Theorems \ref{exampleth1} and \ref{exampleth2}, and Question \ref{question}, we firstly shall give the complete answer to Question \ref{question} in Section 2 (See Theorem \ref{th1} and Theorem \ref{th2}). Besides, In Section 3, we shall also generalize parts $(3)$ and $(4)$ in Theorem \ref{exampleth1} in the similar way as the extension from  parts $(1)$ and $(2)$ in Theorem \ref{exampleth1} to Theorem \ref{exampleth2}, and thus derive several sharp lower and upper bounds for $\K_a-\E_a$ and $\E_a-r^{\thinspace\prime 2}\K_a$ in terms of $\K$, $\E$ and elementary functions.

Let us recall some well-known formulas, which can be found in \cite{AS,AQVV}.
\begin{align}
&\G(x)\G(1-x)=\frac{\pi}{\sin(\pi x)}\quad  (0<x<1), \label{Gamma(1)}\\
&\p(1-x)-\p(x)=\pi \cot(\pi x) \quad  (0<x<1), \label{psi(2)} \\
&\G(x+1)=x\G(x), \quad \p(x+1)=\p(x)+\frac{1}{x} \quad (x>0), \label{GammaPsi} \\
&F(a,b;c;1)=\frac{\G(c)\G(c-a-b)}{\G(c-a)\G(c-b)}\quad (a+b<c), \label{F22}\\
&F(a,b;a+b;x)\sim -\frac{1}{B(a,b)}\log(1-x), \quad  x\rightarrow 1.\label{F33}
\end{align}
Let $\g=\lim_{n\to \infty} (\sum_{k=1}^{\infty}\frac1k-\log n)=0.577215\cdots$ be the Euler-Mascheroni constant. Then the psi function has the following representation (cf. \cite[6.3.16]{AS})
\begin{align}
\p(x)=-\g-\frac{1}{x}+\su \frac{x}{n(n+x)}. \label{psi1}
\end{align}

The following two technical lemmas are useful for proving the monotonicity of functions.

\begin{lemma}\label{LBD}\rm{}
{\rm(\cite[Lemma 5.1]{AQVV})} For $-\infty<a<b<\infty$, let $f$, $g$ : $[a,b]\rightarrow \IR$ be continuous on $[a,b]$, and be differentiable on $(a,b)$. Let $g'(x)\neq 0$ on $(a,b)$. If $f'(x)/g'(x)$ is increasing (decreasing) on $(a,b)$, then so are
\begin{align*}
[f(x)-f(a)]/[g(x)-g(a)] \quad {\rm and} \quad [f(x)-f(b)]/[g(x)-g(b)].
\end{align*}
If $f'(x)/g'(x)$ is strictly monotone, then the monotonicity in the conclusion is also strict.
\end{lemma}

\begin{lemma}\label{series}\rm{}
{\rm({\cite[Lemma 2.1]{Sp1}})} For $n\in\IN_0$, let $r_n$ and $s_n$ be real numbers, and let the power series $R(x)=\s r_nx^n$ and $S(x)=\s s_nx^n$ be convergent for $|x|<1$. If $s_n\geq0$ and not all vanish for $n\in\IN_0$, and if $r_n/s_n$ is strictly increasing (decreasing) in $n\in\IN_0$, then the function $x\mapsto R(x)/S(x)$ is strictly increasing (decreasing, respectively) on $(0,1)$.
\end{lemma}

\bigskip
\section{The answer to question \ref{question}}
In this section, we always assume that $c$ is a fixed constant in $(0,\infty)$, and $a\in(0,c/2]$. Let
\begin{align}
\la_1&=\la_1(a)=\frac{a}{a-c}, \quad \la_2=\la_2(a)=\la_1+1
=\frac{c-2a}{c-a} \label{la1-2}, \\
\la_3&=\la_3(a)=a[\p(c-a)-\p(a)], \quad
\la_4=\la_4(a)=\frac{a(2a-c-1)}{(c-a)(1-a)}, \label{la3-4} \\
\la_5&=\la_5(a)=\frac{a\G(c)[\p(c+1-a)-\p(a)]}
{\G(c)-\G(a)\G(c+1-a)}, \label{la5} \\
\la_6&=\la_6(a)=\la_3+\frac{a}{c-a}=a[\p(c+1-a)-\p(a)]. \label{la6}
\end{align}
For $\la \in \IR$ and $|x|<1$, let
\begin{align}
P_1(\la,c,x)&=\left(\frac2c\right)^{\la} \su \frac{(\frac c2,n)^2}{(c,n)n!}x^n,
\quad \overline{P}_1(c,x)=P_1(0,c,x), \label{P1x} \\
P_2(\la,c,x)&=\left(\frac{2}{c}\right)^{\la}
\su \frac{(1-\frac c2)(\frac c2,n-1)(\frac c2,n)}{(c,n)n!} x^n,
\quad \overline{P}_2(c,x)=P_2(0,c,x), \label{P2x} \\
P_3(\la,c,x)&=\left(\frac 2c\right)^{\la}
\su \frac{(\frac c2,n-1)(\frac c2,n)}{(c,n)(n-1)!}x^n,
\quad \overline{P}_3(c,x)=P_3(0,c,x), \label{P3x} \\
P_4(\la,c,x)&=\left(\frac 2c\right)^{\la-1}
\s \frac{(\frac c2,n)^2}{(n+c)(c,n)n!}x^{n+1},
\quad \overline{P}_4(c,x)=P_4(0,c,x). \label{P4x}
\end{align}

Now we state our main results below.

\begin{theorem}\label{th1}
Let ${\la}^{\ast}_3$ be as in Lemma \ref{la1-6}\,(3) and let $P_1(\la,c,x)$, $P_2(\la,c,x)$, $\overline{P}_1(c,x)$ and $\overline{P}_2(c,x)$ be as in (\ref{P1x})-(\ref{P2x}). For each fixed $c\in(0,\infty)$, $a\in(0,c/2]$, $x\in (0,1)$ and $\la\in\IR$, the functions $f_{1,\,\la}$ and $f_{2,\,\la}$ are given in (\ref{fla12}). Then we have the following conclusions:

$(1)$ $f_{1,\,\la}$ is strictly increasing (decreasing) on $(0,c/2]$ if and only if $\la\leq 0$ ($\la\geq1$ if $c\in(0,1]$ or $\la\geq {\la}^{\ast}_3$ if $c\in(1,\infty)$, respectively), with
\begin{equation*}
f_{1,\,\la}(0^+)
=\begin{cases}
0, &\mbox{~ if ~}\la<1,\\
\log \frac{1}{1-x}, &\mbox{~if ~}\la=1,\\
\infty,&\mbox{~if ~}\la>1,
\end{cases}
\quad f_{1,\,\la}\left(\frac{c}{2}\right)=P_1(\la,c,x).
\end{equation*}
In particular, for $c\in(0,1]$ $(c\in(1,\infty))$, $a\in(0,c/2]$ and $x\in (0,1)$,
\begin{align}\label{IneqTh1}
1+\frac{2a}{c} \overline{P}_1(c,x)
\leq F(a,c-a;c;x) \leq 1+
\min \left\{ \overline{P}_1(c,x), a \log \frac{1}{1-x} \right\},
\end{align}
\begin{align*}
\left( 1+\left(\frac{2a}{c}\right)^{{\la}_{3}^*} \overline{P}_1(c,x)
\leq F(a,c-a;c;x) \leq 1+
\min \left\{ \overline{P}_1(c,x), a \log \frac{1}{1-x} \right\}
,respectively \right).
\end{align*}
The first (second) equality holds if and only if $a=c/2$ ($a=c/2$ or $a\to0$, respectively).

$(2)$ If $c\in(0,2)$, then $f_{2,\,\la}$ is strictly increasing (decreasing) on $(0,c/2]$ if and only if $\la\leq {2}/(c-2)$ ($\la \geq0$, respectively), with
\begin{equation*}
f_{2,\,\la}(0^+)
=\begin{cases}
0, &\mbox{~ if ~}\la<0,\\
$x$,  &\mbox{~if ~}\la=0,\\
\infty,&\mbox{~if ~}\la>0,
\end{cases}
\quad f_{2,\,\la}\left(\frac{c}{2}\right)=P_2(\la,c,x).\\
\end{equation*}
If $c\in[2,\infty)$, then there does not exist $\la \in \IR$ for which $f_{2,\,\la}$ is strictly increasing (or decreasing) on $(0,c/2]$. In particular, for $c\in(0,2)$, $a\in(0,c/2]$ and $x\in (0,1)$,
\begin{align}\label{IneqTh2}
1- \min \left\{x, \left(\frac{2a}{c}\right)^{\frac{2}{c-2}}
\overline{P}_2(c,x) \right\}
\leq F(a-1,c-a;c;x) \leq 1-\overline{P}_2(c,x).
\end{align}
The first (second) equality holds if and only if $a=c/2$ or $a\to0$ ($a=c/2$, respectively).
\end{theorem}

\begin{theorem}\label{th2}
Let ${\la}^{\ast}_3$ be as in Lemma \ref{la1-6}\,(3) and let $P_3(\la,c,x)$, $P_4(\la,c,x)$, $\overline{P}_3(c,x)$ and $\overline{P}_4(c,x)$ as in (\ref{P3x})-(\ref{P4x}). For each fixed $c\in(0,\infty)$ and $a\in(0,c/2]$, $x\in(0,1)$ and $\la\in\IR$, the functions $f_{3,\,\la}$ and $f_{4,\,\la}$ are given in (\ref{fla3}) and (\ref{fla4}), respectively. Then we have the following conclusions:

$(1)$ $f_{3,\,\la}$ is strictly increasing (decreasing) on $(0,c/2]$ if and only if $\la\leq-1$ ($\la\geq1$ if $c \in (0,1]$ or $\la\geq {\la}^{\ast}_3$ if $c\in(1,\infty)$, respectively), with
\begin{equation*}
f_{3,\,\la}(0^+)
=\begin{cases}
0, &\mbox{~ if ~}\la<0,\\
$x$,  &\mbox{~if ~}\la=0,\\
\infty, &\mbox{~if ~}\la>0,
\end{cases}
\quad f_{3,\,\la}\left(\frac{c}{2}\right)=P_3(\la,c,x).\\
\end{equation*}
In particular, for $c\in(0,1]$ $(c\in(1,\infty))$, $a\in(0,c/2]$ and $x\in (0,1)$,
\begin{align}\label{IneqTh3}
\frac{2a}{c} \overline{P}_3(c,x)
\leq F(a,c-a;c;x)-F(a-1,c-a;c;x) \leq \frac{c}{2a} \overline{P}_3(c,x),
\end{align}
\begin{align*}
\left( \left(\frac{2a}{c}\right)^{{\la}^{\ast}_3} \overline{P}_3(c,x)
\leq F(a,c-a;c;x)-F(a-1,c-a;c;x) \leq \frac{c}{2a} \overline{P}_3(c,x)
,respectively \right)
\end{align*}
with equality in each instance if and only if $a=c/2$.

$(2)$ $f_{4,\,\la}$ is strictly increasing (decreasing) on $(0,c/2]$ if and only if $\la\leq 1$ ($\la \geq {\la}^{\ast}_3$, respectively), with
 \begin{equation*}
f_{4,\,\la}(0^+)
=\begin{cases}
0, &\mbox{~ if ~}\la<1,\\
\frac{x}{c},  &\mbox{~if ~}\la=1,\\
\infty,&\mbox{~if ~}\la>1,
\end{cases}
\quad f_{4,\,\la}\left(\frac{c}{2}\right)=P_4(\la,c,x).
\end{equation*}
In particular, for $c\in(0,\infty)$, $a\in(0,c/2]$ and $x\in (0,1)$,
\begin{align}\label{IneqTh4}
\max \left\{ \frac{ax}{c}, \left( \frac{2a}{c} \right)^{{\la}^{\ast}_3}
\overline{P}_4(c,x) \right\}
\leq F(a-1,c-a;c;x)-(1-x)F(a,c-a;c;x)\leq \frac{2a}{c} \overline{P}_4(c,x).
\end{align}
The first (second) equality holds if and only if $a=c/2$ or $a\to0$ ($a=c/2$, respectively).
\end{theorem}
The following three lemmas are required in the proofs of Theorems \ref{th1}--\ref{th2}.

\begin{lemma}\label{la1-6} For each fixed $c\in(0,\infty)$, let $\la_i \,( 0\leq i \leq 6)$ be as in (2.1)-(2.4). Then the following statements hold true:

$(1)$ $\widetilde{\la}_{1}\equiv\mathop{\inf}\limits_{a\in(0,\,c/2]}\la_1(a)=-1$;

$(2)$ $\widetilde{\la}_{2}\equiv\mathop{\inf}\limits_{a\in(0,\,c/2]}\la_2(a)=0$;

$(3)$ $\overline{\la}_{3}\equiv\mathop{\sup}\limits_{a\in(0,\,c/2]}\la_3(a)=
\begin{cases}
1, &\mbox{~ if ~}c\in(0,1],\\
\la^{\ast}_{3}, &\mbox{~if ~}c\in(1,\infty),
\end{cases}$
where $\la^{\ast}_{3}\geq 1$;

$(4)$ $\widetilde{\la}_{4}\equiv\mathop{\inf}\limits_{a\in(0,\,c/2]}\la_4(a)=\begin{cases}
c/(c-2), &\mbox{~ if ~}c\in(0,2),\\
-\infty, &\mbox{~if ~}c\in[2,\infty);
\end{cases}$

$(5)$ $\overline{\la}_{5}\equiv\mathop{\sup}\limits_{a\in(0,\,c/2]}\la_5(a)=\begin{cases}
0, &\mbox{~ if ~}c\in(0,2),\\
\infty, &\mbox{~if ~}c\in[2,\infty);
\end{cases}$

$(6)$ $\overline{\la}_{6}\equiv\mathop{\sup}\limits_{a\in(0,\,c/2]}\la_{6}(a)=\la^{\ast}_{3}$.
\end{lemma}
{\it Proof.}
Since the function $a\mapsto\la_1$ and $a\mapsto\la_2$ are both strictly decreasing on $(0,c/2]$ with ranges $[-1,0)$ and $[0,1)$ respectively, then parts $(1)$ and $(2)$ follows.

For part $(3)$, if $c\in(0,1]$, then from $(1.14)$ we obtain
\begin{align}\label{supla3(1)}
\overline{\la}_{3} \geq \lim_{a\to 0^+ }\la_3(a)=\lim_{a\to 0^+}a\p(c-a)-\lim_{a\to 0^+}a\p(a)=1.
\end{align}
On the other hand, since $\p$ is strictly increasing on $(0,\infty)$, then by (\ref{psi(2)}) we get
\begin{align}\label{supla3(2)}
\overline{\la}_{3} \leq \sup_{a\in(0,\,1/2]} a[\p(1-a)-\p(a)]
=\sup_{a\in(0,\,1/2]} \frac{\pi a}{\tan (\pi a)}=1.
\end{align}
Consequently, it follows from (\ref{supla3(1)}) and (\ref{supla3(2)}) that $\overline{\la}_{3}=1$ if $c\in(0,1]$. Using
the monotonicity of $\p$ on $(0,\infty)$ again,  we obtain $\la^{\ast}_3\geq1$.

For part $(4)$, if $c\in(0,2)$, then $\la_{4}(c/2)=2/(c-2)$, so that $\widetilde{\la}_{4}\leq2/(c-2)$. On the other hand,
\begin{align*}
\la_4(a)=\frac{a(2a-c-1)}{(c-a)(1-a)}\geq \frac{2}{c-2}
\end{align*}
for $a\in(c/2]$ because it is equivalent to the simple inequality $2-a(3-c)\geq0$ for $a\in(0,c/2]$. This yields $\widetilde{\la}_{4}= 2/(c-2)$ if $c\in(0,2)$. If $c\in[2,\infty)$, then from the representation of $\la_4$ one has $\widetilde{\la}_{4}=-\infty$ immediately.

For part $(5)$, if $c\in[2,\infty)$, then
\begin{align*}
\lim_{a\to1^{+}}\la_5(a)
= \G(c) \left[\p(c)-\p(1)\right]
\lim_{a\to1^{+}} \frac{1}{\G(c)-\G(a)\G(c+1-a)}=\infty,
\end{align*}
which shows that $\overline{\la}_5=\infty$.

If $c\in(0,2)$, then by (\ref{psi1}),
\begin{align}\label{la5sup}
\overline{\la}_5\geq
\la_5(0^+)&=\lim_{a\to0^{+}}\frac{a\G(c)[\p(c-a+1)-\p(a)]}
 {\G(c)-\G(a)\G(c-a+1)} \nonumber = \lim_{a\to0^{+}}\frac{a\p(a)}{\G(a)\G(c-a+1)/\G(c)-1} \nonumber \\
&=\lim_{a\to0^{+}} a\p(a)\cdot \lim_{a\to0^{+}} \frac{1}{\G(a)\G(c-a+1)/\G(c)-1}=0.
\end{align}
On the other hand, set $\la_5(a)=p_1(a)/[1-p_2(a)]$, where $p_1(a)=a[\p(c-a+1)-\p(a)]$ and $p_2(a)=\G(a)\G(c-a+1)/\G(c)$. Since $p_2(a)>0$ for $a\in(0,c/2]$, logarithmic differentiation leads to
\begin{align*}
\frac{p'_2(a)}{p_2(a)}=\p(a)-\p(c-a+1)<0
\end{align*}
for $a\in(0,c/2]$, which shows that $p_2$ is strictly decreasing on $a\in(0,c/2]$. Thus $1-p_2(a)\leq 1-p_{2}(1^{-})=0$ for $a\in(0,c/2]$, and therefore $\la_5(a)\leq0$ for $a\in(0,c/2]$. This, in conjunction with (\ref{la5sup}), implies that $\overline{\la}_5=0$.

For part $(6)$, by (\ref{la3-4}), (\ref{la6}) and part (2), we clearly see that
\begin{align*}
\overline{\la}_6
&=\sup_{a\in(0,\,c/2],\,c\in(0,\,\infty)} a\left[ \p(c+1-a)-\p(a) \right] \\
&=\sup_{a\in(0,\,c/2],\,c\in(1,\,\infty)} a\left[ \p(c-a)-\p(a) \right]
={\la}^{\ast}_3.
\end{align*}

\begin{lemma}\label{lemmaanbn}
For each fixed $c\in(0,\infty)$, let $a\in (0,c/2]$ and $n\in\IN_0$, $a_n=\psi(n+a)-\psi(n+c-a)$ and $b_n=\psi(n+a)-\psi(n+1+c-a)$.
Then the sequences $\{a_n\}$ and $\{b_n\}$ are both strictly increasing in $n\in\IN_0$ with $a_{\infty}=\lim_{n\to\infty}a_n=b_{\infty}=\lim_{n\to\infty}b_n=0$.
\end{lemma}

{\it Proof.} By the asymptotic formula for the psi function (cf. \cite[6.3.18]{AS}), $a_{\infty}=b_{\infty}=0$. By (\ref{GammaPsi}), we have
\begin{align*}
a_{n+1}-a_n=\frac{c-2a}{(a+n)(c-a+n)}
\quad \mbox{and} \quad b_{n}=a_{n}-\frac{c}{c-a+n},
\end{align*}
from which the monotonicity properties of $\{a_n\}$ and $\{b_n\}$ follow. \square

\begin{lemma}\label{lemmag1-g4}
For each fixed $c\in(0,\infty)$, let $a\in (0,c/2]$ and $n\in\IN_0$. Define
\begin{equation*}
g_{1,\,n}(a)=(a,n)(c-a,n+1)
\end{equation*}
\begin{equation*}
g_{2,\,n}(a)=(a,n)(c-a,n)
\end{equation*}
\begin{equation*}
g_{3,\,n}(a)=(1-a)(a,n)(c-a,n+1)
\end{equation*}
and
\begin{equation*}
g_{4,\,n}(a)=a(a,n)(c-a,n).
\end{equation*}
Then the sequence $\{g'_{i,\,n}(a)/g_{i,\,n}(a)\}(i=1,2,3,4)$ is strictly increasing in $n\in\IN_0$.
\end{lemma}

{\it Proof.} Let $a_n$ and $b_n$ be as in Lemma \ref{lemmaanbn}. By logarithmic differentiations, we obtain
\begin{equation}\label{g1'/g1}
\frac{g_{1,\,n}'(a)}{g_{1,\,n}(a)}=a_{n+1}-a_0-\frac{1}{n+a},
\end{equation}
\begin{equation}\label{g2'/g2}
\frac{g_{2,\,n}'(a)}{g_{2,\,n}(a)}=a_n-a_0,
\end{equation}
\begin{equation}\label{g3'/g3}
\frac{g_{3,\,n}'(a)}{g_{3,\,n}(a)}=a_{n+1}-a_0-\frac{1}{n+a}+\frac{1}{a-1},
\end{equation}
\begin{equation}\label{g4'/g4}
\frac{g_{4,\,n}'(a)}{g_{4,\,n}(a)}=a_n-a_0+\frac1a.
\end{equation}

Therefore, Lemma \ref{lemmag1-g4} follows from (\ref{g1'/g1})-(\ref{g4'/g4}) and Lemma \ref{lemmaanbn}. \square

\medskip
\noindent {\bf Proof of Theorem \ref{th1}} For part $(1)$, by (\ref{2F1}) and (\ref{GammaPsi}) we obtain
\begin{align}\label{P1cx}
f_{1,\,\la}\left(\frac c2\right)=\left(\frac2c\right)^{\la} \su \frac{(c/2,n)^2}{(c,n)n!}x^n=P_1(\la,c,x)
\end{align}
and
\begin{equation*}
f_{1,\,\la}(a)=\frac{a^{1-\la}\G(c)}{\G(a+1)\G(c-a)}
\su\frac{\G(n+a)\G(n+c-a)}{\G(n+c)n!}x^n,
\end{equation*}
from which it follows that
\begin{align}\label{P1cx0}
f_{1,\,\la}(0^+)&=\lim_{a\to0^+}a^{1-\la}\su \frac{x^n}{n}
=\begin{cases}
0, &\mbox{~ if ~}\la<1,\\
\log \frac{1}{1-x}, &\mbox{~if ~}\la=1,\\
\infty, &\mbox{~if ~}\la>1.
\end{cases}
\end{align}

Logarithmic differentiation yields
\begin{align}\label{f1'}
a\frac{f_{1,\,\la}^{\thinspace\prime}(a)}{f_{1,\,\la}(a)}=F_{1,\,a}(x)-\la,
\end{align}
where
\begin{align}\label{F(1,a,r)}
F_{1,\,a}(x)=a\frac{\frac{\partial}{\partial a} \left[F(a,c-a;c;x)-1\right]}
{F(a,c-a;c;x)-1}=a\frac{\s \frac{g'_{2,n+1}(a)}{(c,n+1)(n+1)!} x^n}
{\s \frac{g_{2,n+1}(a)}{(c,n+1)(n+1)!} x^n},
\end{align}
and $g_{2,\,n}$ is defined in Lemma \ref{lemmag1-g4}. Noting that $g_{2,\,n+1}>0$ for $n\in\IN_0$, and $ g'_{2,\,n+1}/g_{2,\,n+1}$ is strictly increasing in $n\in\IN_0$ by Lemma \ref{lemmag1-g4}. Then by application of Lemma \ref{series}, we derive that $F_{1,\,a}$ is strictly increasing in $x\in(0,1)$. Furthermore, by (\ref{g2'/g2}) and (\ref{F(1,a,r)}), $F_{1,\,a}(0^+)=a g'_{2,\,1}(a)/g_{2,\,1}(a)=(c-2a)/(c-a)=\la_2$. For the limiting value of $F_{1,\,a}(x)$ at $1$, firstly, written $F_{1,\,a}(x)$ as
\begin{align}\label{F1a(x)}
F_{1,\,a}(x)&=\frac{1}{1-1/F(a,c-a;c;x)}F_{2,\,a}(x),
\end{align}
where
\begin{align}\label{F2a}
F_{2,\,a}(x)=a\frac{\frac{\partial}{\partial a}F(a,c-a;c;x)}{F(a,c-a;c;x)}.
\end{align}

Next, for $a\in(0,c/2]$ and $x\in(0,1)$, set
\begin{align*}
F_3(a,x)&=\int_0^{\pi/2}(\sin u)^{2c-2a-1}(\cos u)^{2a-1}(1-x\sin^2 u)^{-a} du,\\
F_4(a,u)&=(\sin u)^{2c-2a-1}(\cos u)^{-1} \log (\sin u)
\end{align*}
and
\begin{align*}
F_5(a,u)&=(\sin u)^{2c-2a-1}(\cos u)^{2a-1}(1-x\sin^2 u)^{-a}\left[2\log(\tan u)+\log\left(1-x\sin^2 u\right)\right].
\end{align*}
Then $F(a,c-a;c;x)=[{2\G(c)}/{\G(a)\G(c-a)}]F_3(a,x)$ by (\ref{int}), so that
\begin{align}
F_{2,\,a}(x)&=\frac{2a\G(c)}{\G(a)\G(c-a)F(a,c-a;c;x)}
 \left\{\frac{\partial F_3}{\partial a}+[\psi(c-a)-\psi(a)]F_3\right\}
 \nonumber \\
&=a[\psi(c-a)-\psi(a)]+\frac{2a\G(c)}{\G(a)\G(c-a)F(a,c-a;c;x)}
\frac{\partial F_3}{\partial a}, \label{F1} \\
\frac{\partial F_3}{\partial a}&=-\int_0^{\pi/2}F_5(a,u)du, \quad
\lim_{x\to1}\frac{\partial F_3}{\partial a}=-2\int_0^{\pi/2}F_4(a,u)du. \label{PF2(a,1)}
\end{align}

It is well known that for $n\in\IN$ (cf. \cite[6.4.1]{AS}),
\begin{align*}
\psi^{(n)}(x)=(-1)^{n+1}\int_0^{\infty}\frac{t^ne^{-xt}}{1-e^{-t}}dt
=-\int_0^1\frac{t^{x-1}}{1-t}(\log t)^ndt.
\end{align*}
Employing (\ref{PF2(a,1)}) and substituting $t=\sin^2u$, we have
\begin{align*}
\lim_{r\to1}\frac{\partial F_3}{\partial a}
&=-2\int_0^{\pi/2}\frac{(\sin u)^{2c-2a-1}\log(\sin u)}{1-\sin^2u}d(\sin u)\\
&=-\frac{1}{2}\int_0^{\pi/2}\frac{(\sin u)^{2(c-a-1)}
\log(\sin^2 u)}{1-\sin^2u}d(\sin^2 u)\\
&=-\frac{1}{2}\int_0^1\frac{t^{(c-a)-1}\log t}{1-t}dt=\frac{1}{2}\psi'(c-a).
\end{align*}
Finally, according to (\ref{F33}), $(\ref{F1a(x)})$ and (\ref{F1}) we get
\begin{align}\label{F2(1)}
F_{1,\,a}(1^-)=F_{2,\,a}(1^-)=a[\psi(c-a)-\psi(a)]=\la_3.
\end{align}

In conclusion, $F_{1,\,a}$ is strictly increasing from $(0,1)$ onto $(\la_2,\la_3)$. Combining (\ref{f1'}) with Lemma \ref{la1-6}\,(1) and (3), we obtain that, for all $c\in(0,\infty)$, $a\in(0,c/2]$ and $x\in(0,1)$,
\begin{equation*}
f_{1,\,\la}^{\thinspace\prime}(a)\geq0 \Longleftrightarrow
\la \leq\inf_{a\in(0,c/2],\,x\in(0,1)}
\{F_{2,\,a}(x)\}=\inf_{a\in(0,c/2]} \{\la_2(a)\}=\widetilde{\la}_2=0,
\end{equation*}
\begin{equation*}
f_{1,\,\la}^{\thinspace\prime}(a)\leq0 \Longleftrightarrow
\la \geq\sup_{a\in(0,c/2],\,x\in(0,1)} \{F_{2,\,a}(x)\}
=\sup_{a\in(0,c/2]} \{\la_3(a)\}=\overline{\la}_3
=\begin{cases}
1,  &\text{if} \,\,\, 0<c\leq 1, \\
{\la}^{*}_3, &\text{if} \,\,\, c>1.
\end{cases}
\end{equation*}
This, together with \eqref{P1cx} and \eqref{P1cx0}, yields the first assertion in part (1). Employing the monotonicity properties and ranges of $f_{1,\,0}$, $f_{1,\,1}$ and $f_{1,\,{\la}^{*}_3}$, we get (\ref{IneqTh1}) immediately.

For part (2), by (\ref{2F1}) and (\ref{GammaPsi}), one can easily obtain
\begin{align}\label{f2(0)11}
f_{2,\,\la}\left(\frac c2\right)&=\left(\frac{2}{c}\right)^{\la}
\su \frac{(1-c/2)(c/2,n-1)(c/2,n)}{(c,n)n!} x^n=P_2(\la,c,x)
\end{align}
and
\begin{align*}
f_{2,\,\la}(a)&=-a^{-\la}\su \frac{(a-1,n)(c-a,n)}{(c,n)n!}x^n\\
&=-a^{-\la}\su\frac{\G(n+a-1)\G(n+c-a)}{\G(a-1)\G(c-a)(c,n)n!}x^n \\
&=(1-a)a^{1-\la}\su\frac{\G(n+a-1)\G(n+c-a)}{\G(a+1)\G(c-a)(c,n)n!}x^n\\
&=(1-a)a^{-\la}\left[\frac{c-a}{c}x
+a\sum_{n=2}^{\infty}\frac{\G(n+a-1)\G(n+c-a)}{\G(a+1)\G(c-a)(c,n)n!}x^n\right],
\end{align*}
from which it follows that
\begin{align}\label{f2(0)1}
f_{2,\,\la}(0^+)&=
\begin{cases}
0, &\mbox{~ if ~}\la<0, \\
x, &\mbox{~ if ~}\la=0,\\
\infty, &\mbox{~ if ~}\la>0.
\end{cases}
\end{align}

Logarithmic differentiation gives
\begin{align}
a\frac{f_{2,\,\la}^{\thinspace\prime}(a)}{f_{2,\,\la}(a)}
&=F_{6,\,a}(x)-\la, \label{f2la'}
\end{align}
where
\begin{align}\label{F6ax1}
F_{6,\,a}(x)=a\frac{\frac{\partial}{\partial a}\su \frac{(1-a)(a,n-1)(c-a,n)}{(c,n)n!} x^n}
{\su \frac{(1-a)(a,n-1)(c-a,n)}{(c,n)n!} x^n}
=a\frac{\s \frac{g_{3,n}'(a)}{(c,n+1)(n+1)!} x^n}
{\s \frac{g_{3,n}(a)}{(c,n+1)(n+1)!} x^n},
\end{align}
and $g_{3,\,n}$ is defined in Lemma \ref{lemmag1-g4}. Since  $g_{3,\,n}>0$ for $n\in\IN_{0}$ and $g'_{3,n}/g_{3,n}$ is strictly increasing in $n\in\IN_{0}$ by Lemma \ref{lemmag1-g4}, then applying Lemma 1.6 we derive that $F_{6,a}$ is strictly increasing in $x\in(0,1)$. Furthermore, by (\ref{g3'/g3}) and (\ref{F6ax1}), we obtain
\begin{align}\label{F6(0)1}
F_{6,\,a}(0^+)&=a\frac{g_{3,\,0}'(a)}{g_{3,\,0}(a)}
=\frac{a(2a-c-1)}{(c-a)(1-a)}=\la_4.
\end{align}

For the limiting value of $F_{6,\,a}(x)$ at 1, firstly, write $F_{6,\,a}(x)$ as
\begin{align}\label{F6ax2}
F_{6,\,a}(x)=a\frac{\frac{\partial}{\partial a}[1-F(a-1,c-a;c;x)]}{1-F(a-1,c-a;c;x)}=\frac{1}{1-1/F(a-1,c-a;c;x)}  F_{7,\,a}(x),
\end{align}
where
\begin{align}\label{F7}
F_{7,\,a}(x)=a\frac{\frac{\partial}{\partial a}F(a-1,c-a;c;x)}{F(a-1,c-a;c;x)}.
\end{align}

Next, for $a\in(0,c/2]$ and $x\in(0,1)$, let
\begin{align*}
F_8(a,x)&=\int_0^{\pi/2}(\sin u)^{2c-2a-1}(\cos u)^{2a-1}(1-x\sin^2 u)^{1-a}du, \\
F_9(a)&=\int_0^{\pi/2}(\sin u)^{2c-2a-1}(\cos u)\log(\sin u)du.
\end{align*}
Then $F(a-1,c-a;c;x)=\left[{2\G(c)}/({\G(a)\G(c-a)})\right]F_8(a,x)$ by (\ref{int}), we obtain
\begin{align*}
\frac{\partial F_8}{\partial a}=-\int_0^{\pi/2}(\sin u)^{2c-2a-1}(\cos u)^{2a-1}
(1-x\sin^2 u)^{1-a}\left[2\log(\tan u)+\log\left(1-x\sin^2 u\right)\right]du
\end{align*}
and
\begin{align}\label{F51}
F_{7,\,a}(x)&=\frac{2a\G(c)}{\G(a)\G(c-a)F(a-1,c-a;c;x)}
\left\{\frac{\partial F_8}{\partial a}+[\p(c-a)-\p(a)]F_8 \right\}\nonumber \\
&=a[\p(c-a)-\p(a)]+\frac{2a\G(c)}{\G(a)\G(c-a)F(a-1,c-a;c;x)}
\frac{\partial F_8}{\partial a}.
\end{align}
It follows from (\ref{F22}) and (\ref{F51}) that
\begin{align}\label{F9a}
F_{7,\,a}(1^-)&=a[\p(c-a)-\p(a)]+\lim_{x\rightarrow 1^-}\frac{2a\G(c)}{\G(a)\G(c-a)F(a-1,c-a;c;1)}
\frac{\partial F_8}{\partial a} \nonumber \\
&=a[\p(c-a)-\p(a)]-4a(c-a)F_9(a).
\end{align}
Using the substitution $t=\sin u$, and integrating by parts, we obtain
\begin{align}\label{F9a1}
F_9(a)&=\int_0^1 t^{2c-2a-1}\log t\, dt=\frac{1}{2(c-a)}\int_0^1\log t \,
d t^{2(c-a)}=-\frac{1}{4(c-a)^2}.
\end{align}
Hence by (\ref{F22}), (\ref{F6ax2}), (\ref{F9a}) and (\ref{F9a1}), we obtain
\begin{align}\label{F6(1-)}
F_{6,\,a}(1^-)&=\frac{1}{1-1/F(a-1,c-a;c;1)}\cdot F_{7,\,a}(1^-) \nonumber \\
&=\frac{a\G(c)[\p(c-a+1)-\p(a)]}{\G(c)-\G(a)\G(c-a+1)}=\la_5.
\end{align}

In conclusion, $F_{6,\,a}$ is strictly increasing from $(0,1)$ onto $(\la_4,\la_5)$. By (\ref{f2la'}) and Lemma \ref{la1-6}\,(3) and (4), we obtain that, for all $c\in(0,\infty)$, $a\in(0,c/2]$ and $x\in(0,1)$,
\begin{equation*}
f_{2,\,\la}^{\thinspace\prime}(a)\geq0 \Longleftrightarrow
\la \leq\inf_{a\in(0,\,c/2],\,x\in(0,1)} \{F_{6,\,a}(x)\}
=\inf_{a\in(0,\,c/2]} \{\la_4(a)\}=\widetilde{\la}_4=
\begin{cases}
\frac{2}{c-2}  &\text{if} \,\,\, 0<c<2, \\
-\infty,  &\text{if} \,\,\, c\geq 2,
\end{cases}
\end{equation*}
\begin{equation*}
f_{2,\,\la}^{\thinspace\prime}(a)\leq0 \Longleftrightarrow
\la \geq \sup_{a\in(0,\,c/2],\,x\in(0,1)} \{F_{6,\,a}(x)\}
=\sup_{a\in(0,\,c/2]} \{\la_5(a)\}=\overline{\la}_5=\begin{cases}
0,  &\text{if} \,\,\, 0<c<2, \\
\infty, &\text{if} \,\,\, c\geq2.
\end{cases}
\end{equation*}
This, together with (\ref{f2(0)11})--(\ref{f2(0)1}), yields the the first assertion in part (2).

Inequality (\ref{IneqTh2}) follows from the  monotonicity properties and ranges of the particular cases $f_{2,\,0}$ and $f_{2,\,{2}/(c-2)}$. The condition of each equality in (\ref{IneqTh2}) is clear. This completes the proof. \square

\medskip
\noindent {\bf Proof of Theorem \ref{th2}} For part $(1)$, by (\ref{2F1}) and (\ref{GammaPsi}), we have
\begin{align}\label{Ka-Ea}
f_{3,\,\la}(a)&=a^{-\la} \s\frac{(a,n)(c-a,n)}{(c,n)n!}x^n
-a^{-\la} \s\frac{(a-1,n)(c-a,n)}{(c,n)n!}x^n \nonumber \\
&=a^{-\la}\su \frac{(a,n-1)(c-a,n)}{(c,n)(n-1)!}x^n,
\end{align}
from which one can easily obtain
\begin{align}
f_{3,\,\la}\left(\frac c2\right)=\left(\frac2c\right)^{\la}
\su \frac{(c/2,n-1)(c/2,n)}{(c,n)(n-1)!}x^n=P_3(\la,c,x),
\end{align}
\begin{align*}
f_{3,\,\la}(a)
&=a^{1-\la}\su \frac{\G(n+a-1)\G(n+c-a)}{\G(1+a)\G(c-a)(c,n)n!}x^n \\
&=a^{-\la}\left[\frac{c-a}{c}x
+a\sum_{n=2}^{\infty}\frac{\G(n+a-1)\G(n+c-a)}{\G(1+a)\G(c-a)(c,n)n!}x^n\right],
\end{align*}
and therefore
\begin{align}\label{limitsOff4}
f_{3,\,\la}(0^+)=
\begin{cases}
0, &\mbox{if} \,\,\, \la<0, \\
x, &\mbox{if} \,\,\, \la=0, \\
\infty, &\mbox{if} \,\,\, \la>0.
\end{cases}
\end{align}

Logarithmic differentiation of $f_{3,\,\la}$ gives
\begin{align}\label{f4'}
a\frac{f_{3,\,\la}^{\thinspace\prime}(a)}{f_{3,\,\la}(a)}=F_{10,\,a}(x)-\la,
\end{align}
where
\begin{align}\label{F10(r)}
F_{10,\,a}(x)=a \frac{ \frac{\partial}{\partial a}\left[F(a,c-a;c;x)-F(a-1,c-a;c;x)\right]}
{F(a,c-a;c;x)-F(a-1,c-a;c;x) }=a \frac{\s \frac{g'_{1,n}(a)}{(c,n+1)n!} x^n}
{\s \frac{g_{1,n}(a)}{(c,n+1)n!} x^n},
\end{align}
and $g_{1,\,n}$ is defined in Lemma \ref{lemmag1-g4}.

Since $g_{1,n}>0$ for $n\in \IN_0$, and $g'_{1,n}/g_{1,n}$ is strictly increasing in $n\in\IN_0$ by Lemma \ref{lemmag1-g4}, then $F_{10,\,a}$ is strictly increasing in $x\in(0,1)$ by application of Lemma \ref{series}. Moreover, $F_{10,\,a}(0^+)=ag'_{1,0}(a)/g_{1,0}(a)=a/(a-c)=\la_1$ by (\ref{F10(r)}), and from (\ref{F22}), (\ref{F33}), (\ref{F2a}), (\ref{F2(1)}), (\ref{F7}) and (\ref{F6(1-)}) we obtain
\begin{align*}
F_{10,\,a}(1^-)&=a\lim_{x\to 1^-} \left\{ \frac{1}{1-\frac{F(a-1,c-a;c;x)}{F(a,c-a;c;x)}}
\Bigg[ \frac{\frac{\partial }{\partial a}F(a,c-a;c;x)}{F(a,c-a;c;x)}
-\frac{\frac{\partial}{\partial a}F(a-1,c-a;c;x)}{F(a,c-a;c;x)} \Bigg]\right\} \\
&=\lim_{x\to 1^-}F_{2,\,a}(x)-\lim_{x\to 1^-}F_{7,\,a}(x)
\lim_{x\to 1^-} \frac{F(a-1,c-a;c;x)}{F(a,c-a;c;x)} \\
&=\lim_{x\to 1^-}F_{2,\,a}(x)=\la_3.
\end{align*}

In conclusion, $F_{10,\,a}$ is strictly increasing from $(0,1)$ onto $(\la_1,\la_3)$. Incorporated with Lemma \ref{la1-6}\,(1) and (3), equation (\ref{f4'}) gives that, for all $c\in(0,\infty)$, $a\in(0,c/2]$ and $x\in(0,1)$,
\begin{equation*}
f_{3,\,\la}^{\thinspace\prime}(a)\geq0 \Longleftrightarrow
\la\leq\inf_{a\in(0,c/2],\,x\in(0,1)} \{F_{10,\,a}(x)\}
=\inf_{a\in(0,c/2]} \{\la_1(a)\} =\widetilde{\la}_1=-1,
\end{equation*}
\begin{equation*}
f_{3,\,\la}^{\thinspace\prime}(a)\leq0 \Longleftrightarrow
\la\geq\sup_{a\in(0,c/2],\,x\in(0,1)} \{F_{10,\,a}(x)\}=\sup_{a\in(0,c/2]}\{\la_3(a)\}
=\overline{\la}_3=\begin{cases}
1,  &\text{if} \,\,\, 0<c\leq 1, \\
{\la}_{3}^*, &\text{if} \,\,\,  c>1.
\end{cases}
\end{equation*}
Therefore, the first assertion in part (1) holds true.
Applying the monotonicity properties of $f_{4,\,-1}$, $f_{4,\,1}$ and $f_{4,\,{\la}_3^*}$, inequality (\ref{IneqTh3}) and its equality cases follow immediately.

\medskip
For part (2), by (\ref{2F1}), we obtain
\begin{align}\label{Ea-xKa}
f_{4,\,\la}(a)&=a^{-\la}\left[  \s\frac{(a-1,n)(c-a,n)}{(c,n)n!}x^n
-(1-x)\s\frac{(a,n)(c-a,n)}{(c,n)n!}x^n \right]  \nonumber \\
&=a^{1-\la}\s \frac{1}{n+c}\frac{(a,n)(c-a,n)}{(c,n)n!}x^{n+1},
\end{align}
from which it follows that
\begin{align*}\label{f4limitsc2}
f_{4,\,\la}\left(\frac{c}{2}\right)=\left(\frac c2\right)^{1-\la}
\s \frac{(c/2,n)^2}{(n+c)(c,n)n!}x^{n+1}=P_4(\la,c,x),
\end{align*}
\begin{align}
f_{4,\,\la}(a)
&=a^{1-\la}\s\frac{\G(n+a)\G(n+c-a)}{(n+c)\G(a)\G(c-a)(c,n)n!}x^{n+1} \nonumber\\
&=a^{1-\la}x \left[\frac1c+\frac{a}{\G(a+1)\G(c-a)}
\su \frac{\G(n+a)\G(n+c-a)}{(n+c)(c,n)n!}x^n \right],
\end{align}
and therefore
\begin{align*}
f_{4,\,\la}(0^+)=\lim_{a\to 0^+}f_{4,\,\la}(a)=
\begin{cases}
0, &\mbox{if}\,\,\,   \la<1,\\
\frac{x}{c}, &\mbox{if}\,\,\,   \la=1,\\
\infty, &\mbox{if}\,\,\,  \la>1.
\end{cases}
\end{align*}

By logarithmic differentiation, we have
\begin{align}\label{f5'}
a\frac{f_{4,\,\la}^{\thinspace\prime}(a)}{f_{4,\,\la}(a)}=F_{11,\,a}(x)-\la,
\end{align}
where
\begin{align}\label{F11}
F_{11,\,a}(x)&=a \frac{\frac{\partial}{\partial a}\left[F(a-1,c-a;c;x)-(1-x)F(a,c-a;c;x)\right]}
{F(a-1,c-a;c;x)-(1-x)F(a,c-a;c;x)}
=a\frac{\s \frac{g'_{4,n}(a)}{(n+c)(c,n)n!} x^n}
{\s \frac{g_{4,n}(a)}{(n+c)(c,n)n!} x^n},
\end{align}
and $g_{4,\,n}$ is defined in Lemma \ref{lemmag1-g4}\,(4).

Since $g_{4,n}>0$ for $n\in \IN_0$, and $g'_{4,n}/g_{4,n}$ is strictly increasing in $n\in\IN_0$ by Lemma \ref{lemmag1-g4}, then $F_{11,\,a}$ is strictly increasing in $x\in(0,1)$ by application of Lemma \ref{series}. Furthermore, by (\ref{F11}) $F_{11,\,a}(0^+)=a({g'_{4,\,0}(a)})/{g_{4,\,0}(a)}=1$, and from from (\ref{F22}), (\ref{F33}), (\ref{F2a}), (\ref{F2(1)}), (\ref{F7}) and (\ref{F6(1-)}) we have
\begin{align*}
F_{11,\,a}(1^-)&=a\lim_{x\to1^-} \Bigg\{ \frac{1}
{1-(1-x)\frac{F(a,c-a;c;x)}{F(a-1,c-a;c;x)}} \\
&\quad \times \left[ \frac{\frac{\partial}
{\partial a}F(a-1,c-a;c;x)}{F(a-1,c-a;c;x)}-(1-x)\frac{\frac{\partial}{\partial a}F(a,c-a;c;x)}{F(a-1,c-a;c;x)}\right]\Bigg\} \\
&=a\lim_{x\to1^-} \frac{\frac{\partial}{\partial a} F(a-1,c-a;c;x)}
{F(a-1,c-a;c;x)}-a\lim_{x\to1^-} \Bigg\{
\frac{(1-x)F(a,c-a;c;x)}{F(a-1,c-a;c;x)}
\frac{\frac{\partial}{\partial a} F(a,c-a;c;x)}{F(a,c-a;c;x)} \Bigg\} \\
&=\lim_{x\to1^-}F_{7,\,a}(x)-\lim_{x\to1^-}\frac{(1-x)F(a,c-a;c;x)}
{F(a-1,c-a;c;x)} \lim_{x\to1^-}F_{2,\,a}(x)=\la_6.
\end{align*}
Therefore, $F_{11,\,a}$ is strictly increasing from $(0,1)$ onto $(1,\la_6)$ by Lemma \ref{series}. Using (\ref{f5'}) and Lemma \ref{la1-6}\,(5), we obtain that, for all $c\in(0,\infty)$, $a\in(0,c/2]$ and $x\in(0,1)$,
\begin{align*}
f_{4,\,\la}^{\thinspace\prime}(a)&\geq0 \Longleftrightarrow \la \leq\inf_{a\in(0,c/2],\,x\in(0,1)} \{F_{11,\,a}(x)\}
=F_{11,\,a}(0^+)=1,\\
f_{4,\,\la}^{\thinspace\prime}(a)&\leq0 \Longleftrightarrow \la\geq\sup_{a\in(0,\,c/2],\,x\in(0,1)} \{F_{11,\,a}(x)\}
=\sup_{a\in(0,\,c/2],\,c\in(0,\,\infty)} \{\la_6(a)\}=\overline{\la}_6={\la}_3^*,
\end{align*}
which yields the monotonicity properties of $f_{4,\,\la}$. The remaining conclusions are obvious.  \square

\bigskip
\section{Monotonicity of generalized elliptic integrals with respect to $a$}

The purpose of this section is to generalize parts (3) and (4) in Theorem \ref{exampleth1}, thus find the analogous extension from parts (1) and (2) in Theorem
\ref{exampleth1} to Theorem \ref{exampleth2}. For $a\in(0,1/2]$, let
\begin{equation}\label{g5a}
g_{5,\,n}(a)\equiv g_{1,\,n}(a)|_{c=1}=(a,n)(1-a,n+1),
\end{equation}
\begin{equation}\label{g6a}
g_{6,\,n}(a)\equiv g_{4,\,n}(a)|_{c=1}=a(a,n)(1-a,n),
\end{equation}
\begin{equation}\label{la7}
\la_7=\la_7(a)=\pi a/\tan(\pi a),
\end{equation}
\begin{equation}\label{la8}
\la_8=\la_8(a)=a[\p(n+a+1)-\p(n+3-a)+\p(1-a)-\p(a)],
\end{equation}
\begin{equation}\label{la9}
\la_9=\la_9(a)=1+a[\p(n+1+a)-\p(n+2-a)+\p(1-a)-\p(a)],
\end{equation}
\begin{equation}\label{la10}
\la_{10}=\la_{10}(a)=\frac{a}{1-a}\frac{\sin(\pi a)+\pi(1-a)\cos(\pi a)
-\pi(1-a)^2 \sum_{k=0}^n (k+1)^{-1}(k!)^{-2} g'_{6,\,k}(a)}
{\sin(\pi a)-\pi(1-a)\sum_{k=0}^n (k+1)^{-1}(k!)^{-2} g_{6,\,k}(a)}.
\end{equation}
For $r\in(0,1)$ and $a\in(0,1/2]$, let
\begin{align}
P_{5,\,n}(a,r)&=\frac{\pi}{2}\sum_{k=0}^n \frac{g_{5,\,k}(a)}{k!(k+1)!}r^{2(k+1)},\quad
\overline{P}_{5,\,n}(r)=P_{5,\,n}\left(\frac12,r\right), \label{P5r}\\
P_{6,\,n}(a,r)&=\frac{\pi}{2}\sum_{k=0}^n
\frac{g_{6,\,k}(a)}{k!(k+1)!}r^{2(k+1)}, \quad
\overline{P}_{6,\,n}(r)=P_{6,\,n}\left(\frac12,r\right), \label{P6r} \\
P_{7,\,n}(r)&=\frac{\pi}{2}\sum_{k=n+1}^{\infty}\frac{1}{k}r^{2(k+1)}
=\begin{cases} \label{P7r}
-\pi r^2\log r^{\thinspace\prime}, &\mbox{if ~} n=0, \\
-\pi r^2\log r^{\thinspace\prime} -\frac{\pi}{2}\sum_{k=1}^{n}\frac{1}{k}r^{2(k+1)}, &\mbox{if ~} n\geq 1,
\end{cases}\\
P_{8,\,n}(r)&=
\begin{cases}\label{P8r}
\frac{\pi}{2}\left(r^2+2r'^2\log r^{\thinspace\prime}\right), &{\rm if ~} n=0, \\
\frac{\pi}{2}\left(r^2+2r'^2\log r^{\thinspace\prime}\right)-\frac{\pi}{2}
\sum_{k=1}^{n}\frac{1}{k(k+1)}r^{2(k+1)}, &{\rm if ~} n\geq1.
\end{cases}
\end{align}

Now we state our main result of this section.
\begin{theorem}\label{th3}
For each $\la\in\IR$, $n\in \IN_0$ and $r\in(0,1)$, let $\overline{P}_{5,\,n}$, $\overline{P}_{6,\,n}$ and $P_{i,\,n}\,(5\leq i \leq8)$ be as in (\ref{P5r})-(\ref{P8r}), let $\widetilde{\la}_8$ and $\overline{\la}_{10}$ be as in Lemma \ref{la7-10}\,(2) and (4), respectively. Define the functions $f_{5,\,\la}$ and $f_{6,\,\la}$ on $(0,1/2]$ by
\begin{align*}
f_{5,\,\la}(a)=\frac{\K_a(r)-\E_a(r)-P_{5,\,n}(a,r)}{a^{\la}} \quad and \quad
f_{6,\,\la}(a)=\frac{\E_a(r)
-r^{\thinspace\prime 2}\K_a(r)-P_{6,\,n}(a,r)}{a^{\la}},
\end{align*}
respectively. Then we have the following conclusions:

$(1)$ $f_{5,\,\la}$ is strictly increasing (decreasing) on $(0,1/2]$ if and only if $\la \leq \widetilde{\la}_8=-1/(2n+3)$ ($\la\geq 1$, respectively), with
\begin{equation*}
f_{5,\lambda}(0^+)=\left\{
  \begin{array}{lll}
    0, & \hbox{$\lambda<1$,} \\
    P_{7,n}(r), &\hbox{$\lambda=1$,} \\
    \infty, & \hbox{$\lambda>1$,}
  \end{array}
\right.
\quad
f_{5,\lambda}\left(\frac{1}{2}\right)=2^{\la}\left[\K(r)-\E(r)-\overline{P}_{5,\,n}(r)\right].
\end{equation*}
In particular, for $n\in \IN_0$, $a\in(0,1/2]$ and $r\in(0,1)$,
\begin{align}\label{In8}
2a \left[\K(r)-\E(r)-\overline{P}_{5,\,n}(r)\right]
&\leq \K_a(r)-\E_a(r)-P_{5,\,n}(a,r) \nonumber\\
\leq \min &\left\{ aP_{7,\,n}(r),\, (2a)^{-\frac{1}{2n+3}}
\left[ \K(r)-\E(r)-\overline{P}_{5,\,n}(r)\right] \right\},
\end{align}
with equality in each instance if and only if $a=1/2$ or $a\to 0$.

$(2)$ $f_{6,\,\la}$ is strictly increasing (decreasing) on $(0,1/2]$ if and only if $\la\leq 1$ ($\la\geq \overline{\la}_{10}$, respectively), with
\begin{equation*}
f_{6,\,\la}(0^+)=\left\{
  \begin{array}{lll}
    0, & \hbox{$\lambda<2$,} \\
    P_{8,\,n}(r), &\hbox{$\lambda=2$,} \\
    \infty, & \hbox{$\lambda>2$,}
  \end{array}
\right.
\quad
f_{6,\,\la}\left(\frac{1}{2}\right)=2^{\la}
\left[\E(r)-r^{\thinspace\prime2}\K(r)-\overline{P}_{6,\,n}(r)\right].
\end{equation*}
In particular, for $n\in \IN_0$, $a\in(0,1/2]$ and $r\in(0,1)$,
\begin{align}\label{In9}
(2a)^{\overline{\la}_{10}} \left[\E(r)-r^{\thinspace\prime2}\K(r)-\overline{P}_{6,\,n}(r)\right]
\leq \E_a(r)&-r^{\thinspace\prime2}\K_a(r)-P_{6,\,n}(a,r)  \nonumber \\
\leq &2a\left[\E(r)-r^{\thinspace\prime2}\K(r)-\overline{P}_{6,\,n}(r)\right],
\end{align}
with equality in each case if and only if $a=1/2$.
\end{theorem}

The proof of Theorem \ref{th3} requires some properties of $\la_i(a)$ ($7\leq i\leq10$), which are given in the following Lemma \ref{la7-10}.

\begin{lemma}\label{la7-10}
For $n\in \IN_0$ and $a\in(0,1/2]$, let $\la_i \,(7\leq i\leq 10)$ be as in (\ref{la7})-(\ref{la10}). Then we have the following conclusions:

$(1)$ $\overline{\lambda}_{7}=\sup\limits_{a\in(0,1/2]}\{\lambda_{7}(a)\}=1$;

$(2)$ $\widetilde{\lambda}_{8}=\inf\limits_{a\in(0,1/2]}\{\lambda_{8}(a)\}=-\frac{1}{2n+3}$;

$(3)$ $\widetilde{\lambda}_{9}=\inf\limits_{a\in(0,1/2]}\{\lambda_{9}(a)\}=1$;

$(4)$ Let $\overline{\lambda}_{10}=\sup\limits_{a\in(0,1/2]}\{\lambda_{10}(a)\}$ for $n\in \IN_0$. Then $\overline{\lambda}_{10}=2$ if $n=0$. That is, for each $a\in(0,1/2]$, if we let
\begin{align*}
z(a)\equiv \la_{10}(a)|_{n=0}=\frac{a}{1-a}
\frac{\sin(\pi a)+\pi(1-a)\cos(\pi a)-\pi(1-a)^2}{\sin(\pi a)-\pi a(1-a)},
\end{align*}
then $\sup_{a\in(0,1/2]}\{z(a)\}=2$.
\end{lemma}
\proof
Since the function $x\mapsto x/\tan{x}$ is strictly decreasing from $(0,\pi/2)$ onto $(0,1)$, then part (1) follows.

For part (2), by (\ref{GammaPsi}), $\la_8(1/2)=-1/(2n+3)$, so that
\begin{align}\label{la4(1)}
\widetilde{\la}_8=\inf_{a\in(0,1/2]}\{\la_8(a)\}\leq -\frac{1}{2n+3}.
\end{align}
On the other hand, it is easy to verify that the function $a\mapsto a/(n+2-a)$ is strictly increasing from $(0,1/2]$ onto $(0,1/(2n+3)]$ for each fixed $n\in \IN_0$. Combining (\ref{GammaPsi}) and Lemma \ref{lemmaanbn} gives
\begin{align}\label{la4(2)}
\la_8(a)&=a\left[ \p(n+a+1)-\p(n+2-a)+\p(1-a)-\p(a)\right]-\frac{a}{n+2-a} \nonumber\\
&\geq \frac{1-2a}{1-a}-\frac{a}{n+2-a}\geq -\frac{1}{2n+3}.
\end{align}
Consequently, it follows from (\ref{la4(1)}) and (\ref{la4(2)}) that $\widetilde{\la}_8=-1/(2n+3)$.

For part (3), similarly, one can easily obtain $\la_9(1/2)=1$, so that
\begin{align}\label{la5(1)}
\widetilde{\la}_9=\inf_{a\in(0,1/2]}\{\la_9(a)\}\leq 1.
\end{align}
On the other hand,  since the function $a\mapsto a/(1-a)$ is strictly increasing from $(0,1/2]$ onto $(0,1]$, combining (\ref{GammaPsi}) and Lemma \ref{lemmaanbn} we obtain
\begin{align}\label{la5(2)}
\la_9(a)&=1+a \left[ \p(n+1+a)-\p(n+2-a)+\p(1-a)-\p(a) \right] \nonumber\\
&\geq 1+a \left[ \p(1+a)-\p(2-a)+\p(1-a)-\p(a) \right] \nonumber\\
&=2-\frac{a}{1-a}\geq 1.
\end{align}
Consequently, it follows from (\ref{la5(1)}) and (\ref{la5(2)}) that $\widetilde{\la}_9=1$ for $n\in \IN_0$.

For part (4), let $z_1(a)=a\left[\sin(\pi a)+\pi(1-a)\cos(\pi a)-\pi(1-a)^2\right]$ and $z_2(a)=(1-a)\left[\sin(\pi a)-\pi a(1-a)\right]$. Then $z(a)=z_{1}(a)/z_{2}(a)$. Utilizing the following series expansions
\begin{align*}
\sin x=\s (-1)^n \frac{x^{2n+1}}{(2n+1)!} \quad  {\rm and} \quad
\cos x=\s (-1)^n \frac{x^{2n}}{(2n)!} \quad  (x\in \IR),
\end{align*}
we derive that  $\sin(\pi a)>\pi a-\pi^3 a^3/6$ and $\cos(\pi a)<1-\pi^2 a^2/2+\pi^4a^4/24$ for $a\in(0,1/2]$,
\begin{equation*}
\lim_{a\mapsto 0^+}z(a)=\lim_{a\mapsto 0^+}\frac{a}{1-a}\left[\frac{\pi a+\pi(1-a)(1-\pi^2a^2/2)-\pi(1-a)^2+o(a^3)}{\pi a-\pi a+\pi a^2+o(a^2)}\right]=2
\end{equation*}
and therefore
\begin{equation*}
z_{2}(a)\geq \sin(\pi a)-\pi a(1-a)>\pi a-\frac{\pi^3}{6}a^3-\pi a(1-a)=\pi a^2(1-\pi^2a/6)>0
\end{equation*}
for $a\in(0,1/2]$.

Following it suffices to show that $z_1(a)<2z_2(a)$ for each $a\in(0,1/2]$, which is equivalent to
$z_3(a)\equiv\pi a(1-a)\cos(\pi a)-(2-3a)\sin(\pi a)+\pi a(1-a)^2<0$
for $a\in(0,1/2]$. Noting that
\begin{align*}
z_3(a)&<\pi a(1-a)(1-\pi^2 a^2/2+\pi^4a^4/24)-(2-3a)(\pi a-\pi^3 a^3/6)+\pi a(1-a)^2 \nonumber \\
&=-\frac{\pi a^3}{24} \left[4\pi^2-24-\pi^4a^2(1-a)\right]<0
\end{align*}
for $a\in(0,1/2]$. This yields the assertion of part (4).
\endproof

\medskip
\noindent {\bf Proof of Theorem \ref{th3}} (1) If $n=0$, then  Theorem \ref{th3}\,(1) has been proved in \cite[Theorem 1.2\,(3)]{QMB}. Now we suppose that $n\in\IN_0$, let $g_{5,\,n}(a)$ and $P_{5,\,n}(a,r)$ be as in (\ref{g5a}) and (\ref{P5r}), respectively.  Let $h_1(a)=\K_a(r)-\E_a(r)-P_{5,\,n}(a,r)$, then by (\ref{Ka}) and (\ref{Ea}),
\begin{align}\label{g3(a)}
h_1(a)&=\frac{\pi}{2}\sum_{k=n+1}^{\infty}\frac{g_{5,\,k}(a)}{k!(k+1)!}r^{2k+2}
=\frac{\pi}{2}\sum_{k=0}^{\infty}\frac{g_{5,\,k+n+1}(a)}
{(k+n+1)!(k+n+2)!}r^{2(k+n+2)},
\end{align}
from which it follows that
\begin{align*}
f_{5,\,\la}(a)&=a^{-\la}h_1(a)
=\frac{\pi}{2}a^{-\la}\sum_{k=0}^{\infty}\frac{g_{5,\,k+n+1}(a)}
{(k+n+1)!(k+n+2)!}r^{2(k+n+2)} \nonumber\\
&=\frac{\pi}{2}a^{1-\la}\sum_{k=0}^{\infty}
\frac{\G(k+n+1+a)\G(k+n+3-a)}{\G(a+1)\G(1-a)(k+n+1)!(k+n+2)!}r^{2(k+n+2)},
\end{align*}
and therefore
\begin{equation}\label{f5limits}
f_{5,\,\la}(0^+)=\left\{
\begin{array}{ll}
    0 & \hbox{if $\la<1$,} \\
    P_{7,\,n}(r), & \hbox{if $\la=1$,} \\
    \infty, & \hbox{if $\la>1$,}
  \end{array}
\right.
\quad f_{5,\,\la}\left(\frac{1}{2}\right)=2^{\la}\left[\K(r)-\E(r)-\overline{P}_{5,\,n}(r)\right].
\end{equation}

Logarithmic differentiation of $f_{5,\,\la}$ leads to
\begin{align}\label{F5}
a\frac{f'_{5,\,\la}(a)}{f_{5,\,\la}(a)}
=a\frac{h_{1}'(a)}{h_{1}(a)}-\la=G_{1,a}(r)-\la,
\end{align}
where
\begin{align}\label{G1(ar)}
G_{1,a}(r)&=a\left\{ \sum_{k=0}^{\infty}
\frac{g_{5,\,k+n+1}'(a)}{(k+n+1)!(k+n+2)!}r^{2k} \right\}
\left\{ \sum_{k=0}^{\infty}\frac{g_{5,\,k+n+1}(a)}
{(k+n+1)!(k+n+2)!}r^{2k} \right\}^{-1} \nonumber \nonumber \\
&=a \left\{ \sum_{k=0}^{\infty} C_k r^{2k} \right\}
\left\{ \sum_{k=0}^{\infty} D_k r^{2k} \right\}^{-1},
\end{align}
\begin{align*}
C_k=\frac{g_{5,\,k+n+1}'(a)}{(k+n+1)!(k+n+2)!},\quad
D_k=\frac{g_{5,\,k+n+1}(a)}{(k+n+1)!(k+n+2)!}.
\end{align*}
Clearly, $D_k>0$ for $k\in \IN_{0}$, and $C_k/D_k$ is strictly increasing in $k\in \IN_0$ by Lemma \ref{lemmag1-g4}.  Applying Lemma \ref{series}, $G_{1,a}$ is strictly increasing in $r\in(0,1)$. Moreover, by (\ref{g5a}) and (\ref{G1(ar)}),
\begin{equation*}
G_{1,a}(0^+)=a\frac{C_0}{D_0}=a\frac{g'_{5,\,n+1}(a)}{g_{5,\,n+1}(a)}=\la_8,
\end{equation*}
and from (\ref{Ka}), (\ref{Ea}) and the proof of \cite[Theorem 1.1]{QMB} one has
\begin{align*}
G_{1,a}(1^-)&=\lim_{r\to1^-}\frac{a}{1-\E_a/\K_a-P_{5,\,n}(a,r)/\K_a}
\left[ \frac{1}{\K_a}\frac{\partial \K_a}{\partial a}
-\frac{1}{\K_a}\frac{\partial \E_a}{\partial a}
-\frac{1}{\K_a}\frac{\partial P_{5,\,n}(a,r)}{\partial a}\right] \nonumber \\
&=\lim_{r\to1^-}\frac{a}{\K_a}\frac{\partial \K_a}{\partial a}=\la_7.
\end{align*}
Combining with
(\ref{F5}), and Lemma \ref{la7-10}\,(1) and (2), we obtain that, for all $a\in(0,1/2]$ and $r\in(0,1)$,
\begin{equation*}
f_{5,\,\la}^{\thinspace\prime}(a) \geq0 \Longleftrightarrow \la\leq\inf_{a\in(0,1/2],\,r\in(0,1)} \{G_{1,a}(r)\}
=\inf_{a\in(0,1/2]}\{\la_8(a)\}=\widetilde{\la}_8=-\frac{1}{2n+3},
\end{equation*}
\begin{equation*}
f_{5,\,\la}^{\thinspace\prime}(a) \leq0 \Longleftrightarrow \la\geq\sup_{a\in(0,1/2],\,r\in(0,1)}\{G_{1,a}(r)\}
=\sup_{a\in(0,1/2]}\{\la_7(a)\}=\overline{\la}_7=1.
\end{equation*}
This, together with \eqref{f5limits}, yields the first assertion of part (1).
The double inequality (\ref{In8}) and its equality case are clear.

(2) For $n\in\IN_0$, let $g_{6,\,n}(a)$ and $P_{6,\,n}(a,r)$ be as in (\ref{g6a}) and (\ref{P6r}), respectively. Define $h_2(a)=\E_a(r)-r^{\thinspace\prime2}\K_a(r)-P_{6,\,n}(a,r)$, then by (\ref{Ka}) and (\ref{Ea}) one has
\begin{align}\label{g4a}
h_2(a)&=\frac{\pi}{2}\sum_{k=n+1}^{\infty}\frac{g_{6,\,k}(a)}{k!(k+1)!}r^{2(k+1)}
=\frac{\pi}{2}\sum_{k=0}^{\infty}\frac{g_{6,\,k+n+1}(a)}
{(k+n+1)!(k+n+2)!}r^{2(k+n+2)},
\end{align}
from which it follows that
\begin{align}\label{f6a}
f_{6,\,\la}(a)&=a^{-\la}h_2(a)
=\frac{\pi}{2}a^{-\la}\sum_{k=0}^{\infty}\frac{g_{6,\,k+n+1}(a)}
{(k+n+1)!(k+n+2)!}r^{2(k+n+2)} \nonumber\\
&=\frac{\pi}{2}a^{2-\la}\sum_{k=0}^{\infty}
\frac{\G(k+n+1+a)\G(k+n+2-a)}{\G(a+1)\G(1-a)(k+n+1)!(k+n+2)!}r^{2(k+n+2)},
\end{align}
and therefore
\begin{equation}\label{f6limits1}
f_{6,\,\la}(0^+)=\left\{
\begin{array}{ll}
    0 & \hbox{if $\la<2$,} \\
    \frac{\pi}{2}\sum_{k=n+1}^{\infty}\frac{1}{k(k+1)}r^{2k+2}=P_{8,\,n}(r), & \hbox{if $\la=2$,} \\
    \infty, & \hbox{if $\la>2$,}
  \end{array}
\right.
\end{equation}
\begin{equation}\label{f6limits2}
f_{6,\,\la}\left(\frac{1}{2}\right)=2^{\la}
\left[\E(r)-r^{\thinspace\prime2}\K(r)-\overline{P}_{6,\,n}(r)\right].
\end{equation}
It is easy to verify that for $x\in(0,1)$,
\begin{align*}
\sum_{k=1}^{\infty}\frac{1}{k(k+1)}x^{k+1}=x-(1-x)\log\frac{1}{1-x},
\end{align*}
and hence $P_{8,\,n}(r)=\frac{\pi}{2}\left(r^2+2r^{\thinspace\prime2}\log r^{\thinspace\prime}\right)$ if $n=0$, and
\begin{align*}
P_{8,\,n}(r)=
\frac{\pi}{2}\left[r^2+2r^{\thinspace\prime2}\log r^{\thinspace\prime}
-\sum_{k=1}^{n}\frac{1}{k(k+1)}r^{2(k+1)}\right] \mbox{ ~if ~} n\geq1.
\end{align*}

Logarithmic differentiation of $f_{6,\,\la}$ leads to
\begin{equation}\label{F6}
a\frac{f'_{6,\,\la}}{f_{6,\,\la}}=a\frac{h_{2}'(a)}{h_{2}(a)}-\la=G_{2,a}(r)-\la,
\end{equation}
where
\begin{align}\label{G2ar}
G_{2,a}(r)&=a\left\{ \sum_{k=0}^{\infty}\frac{g_{6,\,k+n+1}'(a)}
{(k+n+1)!(k+n+2)!}r^{2k} \right\} \left\{ \sum_{k=0}^{\infty}
\frac{g_{6,\,k+n+1}(a)}{(k+n+1)!(k+n+2)!}r^{2k} \right\}^{-1} \nonumber\\
&=a \left\{ \sum_{k=0}^{\infty} E_k r^{2k}\right\}
\left\{ \sum_{k=0}^{\infty} F_k r^{2k} \right\}^{-1},
\end{align}
\begin{align*}
E_k=\frac{g_{6,\,k+n+1}'(a)}{(k+n+1)!(k+n+2)!}, \quad
F_k=\frac{g_{6,\,k+n+1}(a)}{(k+n+1)!(k+n+2)!}.
\end{align*}
Clearly, $F_k>0$ for $k\in \IN_0$, and the sequence  $\{E_k/F_k\}$ is strictly increasing in $k\in \IN_0$ by Lemma \ref{lemmag1-g4}, so that the function $r\mapsto G_{2,a}$ is strictly increasing in $(0,1)$ by application of Lemma \ref{series}. Moreover, by (\ref{g6a}) and (\ref{G2ar}),
\begin{align}\label{F6(0)}
G_{2,a}(0^+)&=a\frac{E_0}{F_0}=a \frac{g_{6,\,n+1}'(a)}{g_{6,\,n+1}(a)}=\la_9(a).
\end{align}
and from (\ref{Ea}) and \cite[equation (4.19)]{QMB} we obtain
\begin{align}
G_{2,a}(1^-)&=\lim_{r\to1^-} \Bigg\{ \frac{a[\E_a(r)-\pi/2]}
{\E_a(r)-r^{\thinspace\prime 2}\K_a(r)-P_{6,\,n}(a,r)} \nonumber\\
&\quad \times
\Bigg[ \frac{\frac{\partial}{\partial a}\left[\pi/2-\E_a(r)\right]}{\pi/2-\E_a(r)}
+\frac{\frac{\partial }{\partial a}
\left[r^{\thinspace\prime 2}\K_a(r)\right] }{\pi/2-\E_a(r)}
+\frac{\frac{\partial }{\partial a} P_{6,\,n}(a,r)}{\pi/2-\E_a(r)} \Bigg]\Bigg\}
\nonumber \\
&=\frac{a[\sin(\pi a)-\pi(1-a)]}
{\sin(\pi a)-\pi(1-a)\sum_{k=0}^n (k+1)^{-1}(k!)^{-2} g_{6,\,k}(a)} \nonumber \\
&\quad \times
\left\{ -\frac{\sin(\pi a)+\pi(1-a)\cos(\pi a)}{(1-a)[\pi(1-a)-\sin(\pi a)]}
+\frac{\pi(1-a)\sum_{k=0}^n (k+1)^{-1}(k!)^{-2} g'_{6,\,k}(a)}
{\pi(1-a)-\sin(\pi a)}\right\} \nonumber  \\
&=\frac{a}{\sin(\pi a)-\pi(1-a)
\sum_{k=0}^n (k+1)^{-1}(k!)^{-2} g_{6,\,k}(a) } \nonumber \\
&\quad \times  \left\{ \frac{\sin(\pi a)+\pi(1-a)\cos(\pi a)}{1-a}
-\pi(1-a)\sum_{k=0}^n (k+1)^{-1}(k!)^{-2} g'_{6,\,k}(a)\right\} \nonumber \\
&=\frac{a}{1-a} \frac{\sin(\pi a)+\pi(1-a)\cos(\pi a)
-\pi(1-a)^2 \sum_{k=0}^n (k+1)^{-1}(k!)^{-2} g'_{6,\,k}(a)}
{\sin(\pi a)-\pi(1-a)\sum_{k=0}^n (k+1)^{-1}(k!)^{-2} g_{6,\,k}(a)} \nonumber \\
&=\la_{10}. \label{F6(1)}
\end{align}
Combining with (\ref{F6}), and Lemma \ref{la7-10}\,(3) and (4), we conclude that for all $a\in(0,1/2]$ and $r\in(0,1)$,
\begin{equation*}
f_{6,\,\la}^{\thinspace\prime}(a)\geq0\Longleftrightarrow \la\leq\inf_{a\in(0,1/2],\,r\in(0,1)}\{G_{2,a}(r)\}
=\inf_{a\in(0,1/2]}\{\la_9(a)\}=\widetilde{\la}_9=1,
\end{equation*}
\begin{equation*}
f_{6,\,\la}^{\thinspace\prime}(a) \leq0\Longleftrightarrow \la\geq\sup_{a\in(0,1/2],\,r\in(0,1)}\{G_{2,a}(r)\}
=\sup_{a\in(0,1/2]}\{\la_{10}(a)\}=\overline{\la}_{10}.
\end{equation*}
This, together with \eqref{f6limits1} and  \eqref{f6limits2}, yields the first assertion of part (2).
Inequality (\ref{In9}) and its equality case are clear. This completes the proof.    \square

\medskip
Taking $n=0$ in  Theorem \ref{th3}\,(2) and using Lemma \ref{la7-10}\,(4), the following Corollary \ref{corollary}
can be obtained immediately.

\begin{corollary}\label{corollary}
For each $\la\in\IR$ and $r\in(0,1)$, define the function $f_{7,\,\la}$ on $(0,1/2]$ by
\begin{align*}
f_{7,\,\la}(a)=\frac{\E_a(r)-r^{\thinspace\prime 2}\K_a(r)-\pi ar^2/2}{a^{\la}}.
\end{align*}
Then $f_{7,\,\la}$ is strictly decreasing on $(0,1/2]$ if and only if $\la \geq 2$. In particular, for all $a\in(0,1/2]$ and $r\in(0,1)$,
\begin{align*}
a^2 \left[4(\E-r^{\thinspace\prime 2}\K)-\pi r^2\right] \leq \E_a(r)-r^{\thinspace\prime 2}\K_a(r)-\frac{\pi}{2}ar^2
\leq \frac{\pi a^2}{2}\left(r^2+2r^{\thinspace\prime 2}
\log r^{\thinspace\prime}\right).
\end{align*}
The first (second) equality holds if and only if $a=1/2$ ($a\to 0$, respectively).
\end{corollary}

\begin{remark}\rm{}
It is not easy to calculate the exact values of ${\la}_{3}^{*}$ and $\overline{\la}_{10}\,(n\in \IN)$ appeared in Lemma \ref{la1-6}\,(3) and Lemma \ref{la7-10}\,(4), respectively. Here it is left as an open problem for the readers.
\end{remark}

\begin{open problem}
$(1)$ What is the value of ${\la}_{3}^{*}$ in Lemma \ref{la1-6}\,(2)?

$(2)$ What is the value of $\overline{\la}_{10}$  in Lemma \ref{la7-10}\,(4) for $n\in \IN$?
\end{open problem}

\section*{Conflict of interest}

The authors declare that they have no conflict of interest.

\end{document}